# BANDWIDTH CHOICE FOR NONPARAMETRIC CLASSIFICATION

By Peter Hall and Kee-Hoon Kang

*Australian National University and Australian National University and Hankuk University of Foreign Studies*

It is shown that, for kernel-based classification with univariate distributions and two populations, optimal bandwidth choice has a dichotomous character. If the two densities cross at just one point, where their curvatures have the same signs, then minimum Bayes risk is achieved using bandwidths which are an order of magnitude larger than those which minimize pointwise estimation error. On the other hand, if the curvature signs are different, or if there are multiple crossing points, then bandwidths of conventional size are generally appropriate. The range of different modes of behavior is narrower in multivariate settings. There, the optimal size of bandwidth is generally the same as that which is appropriate for pointwise density estimation. These properties motivate empirical rules for bandwidth choice.

**1. Introduction.**

1.1. *Motivation and main results.* A common approach to nonparametric classification based on data from training samples is to construct nonparametric estimators of population densities and substitute them for the true densities in a theoretically optimal algorithm for minimizing Bayes risk. Not only is this approach intuitively appealing and operationally straightforward, it is optimal in a minimax sense, as argued by Marron (1983). However, it is unclear how one might select a bandwidth that minimizes risk. In particular, we might ask from a theoretical viewpoint what relationship exists between the sizes of bandwidth that are appropriate for pointwise density estimation and for optimal classification. And even if we understand this connection, and have a theoretically optimal formula for bandwidth, how might we go about constructing empirical approximations to it?

In this note we briefly summarize how bandwidth choice influences classification error, and suggest ways of choosing bandwidth to minimize that error. In particular, we show that when only two populations are involved, when the populations are univariate, and when the densities intersect at a single point, the following dichotomous result arises. If the density curvatures are of different signs at the crossing point, then minimum Bayes risk is achieved using bandwidths that are of the same sizes as those which minimize pointwise estimation error. On the other hand, if the curvatures are of the same sign, then quite different bandwidth sizes, in fact, similar to those that would be employed if the kernel was of fourth (rather than second) order, are appropriate. Furthermore, if there is more than one crossing point, then, generally speaking, the first of these two sizes of bandwidth applies.

Ironically, the problem actually becomes simpler in more complex settings, where the classification problem involves multivariate data. There, it is generally the case that the optimal size of bandwidth (in the sense of minimizing Bayes risk) is the same as that which would be used if we were constructing pointwise density estimators.

The problem of empirical bandwidth choice suffers from unexpected difficulties. It might reasonably be thought that leave-one-out methods, which have been so successful in related problems of nonparametric inference [see, e.g., Hall (1983), Stone (1984), Härdle and Kelly (1987) and Györfi, Kohler, Krzyżak and Walk (2002)], would perform well in this setting. For example, one could compute the estimate of classification error when a given datum $X$ was omitted from the sample, evaluate the estimate at $X$, and then average over all values of $X$ in order to obtain an estimate of classification error that could be minimized with respect to bandwidth. However, we shall show that this generally gives poor performance. The reason is that it depends on properties of density estimators at the relatively small number of places where the true densities cross, and the leave-one-out approach described above does not give consistent estimates of error at individual points such as $x$; it is necessary to average over a continuum of points in the neighborhood of $x$. The extra degree of smoothing required by this step complicates inference, with the result that alternative approaches are relatively attractive.

1.2. *Relationship to literature.* The extensive literature on this topic includes results which, at first sight, might appear to be contradictory. For example, it is known that, while there exists a class of universally consistent classifiers [see, e.g., Lugosi and Nobel (1996)], the convergence rate of any classifier can be arbitrarily slow [Devroye, Györfi and Lugosi (1996), Chapter 7 and Yang (1999a)]. Indeed, arbitrarily slow rates can apply even for smooth densities [Devroye (1982)]. Moreover, while for large classes of densities (e.g., monotone ones) the rate of convergence of the risk for classification is strictly faster than that for estimation, the two problems are, in fact, of the same difficulty in a well-defined sense [Yang (1999a)]. Also, although the risk of members of a popular class of classifiers converges to its asymptotic limit at rate $n^{-2}$, where $n$ denotes sample size [Cover (1968)], that for classifiers based on empirical forms of Bayes risk converges no more quickly than $n^{-1}$, even in parametric settings [e.g., Kharin and Ducinskas (1979)]. If Bayes risk-based classifiers use kernel estimators, or related nonparametric methods based on places where densities cross, then they converge at slower rates than $n^{-1}$, which are nevertheless minimax-optimal [e.g., Marron (1983) and Mammen and Tsybakov (1999)].

Such contrasts, particularly those between results of Lugosi and Nobel (1996) and Devroye, Györfi and Lugosi (1996), or among the convergence-rate results noted by Yang (1999a), are particularly engaging, but, of course, do not amount to contradictions. Differences among minimax results can be accommodated by noting that the classes over which the "max" part of "minimax" is taken are not identical. There is no real conflict between the results of Cover (1968) and those for Bayes risk-based methods, since the limiting risk of the nearest-neighbor methods treated by Cover is (except in degenerate cases) strictly greater than the Bayes risk, and so the fast convergence rate does not imply good performance.

Work in the present paper relates to kernel-based methods for classification, which date from contributions of Fix and Hodges (1951). It is less closely connected to classification problems involving very high-dimensional data; for the latter setting, see, for example, Breiman (1998, 2001), Schapire, Freund, Bartlett and Lee (1998), Friedman, Hastie and Tibshirani (2000), Kim and Loh (2001), Dudoit, Fridlyand and Speed (2002) and Jiang (2002). Although there is some evidence that multiplicative bias/variance decompositions play an important role in such contexts, considerable interest still resides in additive decompositions of the type addressed in the results we shall discuss. For example, in a wide-ranging contribution to classification problems for multivariate (and, in particular, high-dimensional) data, Friedman [(1997), Section 11] draws particular attention to the role of additive decompositions in classification problems.

In addition to the work discussed above, there is an extensive literature on nonparametric methods for classification, much of it based on using an empirical version of the Bayes-optimal rule. Fukunaga and Hummels (1987) and Psaltis, Snapp and Venkatesh (1994) extend Cover's (1968) work to $d$ dimensions, where the classification error of nearest-neighbor methods converges at rate $n^{-2/d}$. Efron (1983) and Efron and Tibshirani (1997) discuss the performance of bootstrap-based estimators of error rate for general classification methods. Chanda and Ruymgaart (1989) address kernel-based classification rules when the two distributions differ only in location, and where tails decrease exponentially fast or in a regularly varying manner. See also Kharin (1983), who gives related results in multivariate settings, and Devroye, Györfi and Lugosi [(1996), Theorem 6.6], who provide an elegant upper bound. Krzyżak (1991) derives bounds on Bayes probability of error for kernel-based classification rules; Lapko (1993) gives a book-length account, in Russian, of nonparametric classification, including techniques based on nonparametric density estimation; Pawlak (1993) proposes kernel-based classification rules for use with incomplete data; Lugosi and Pawlak (1994) describe properties of a posterior-probability estimator of classification error for nonparametric classifiers; Ancukiewicz (1998) introduces class-based classification rules founded on nonparametric density estimators; Yang (1999b) studies nonparametric estimation of conditional probability for classification; Baek and Sung (2000) introduce a nearest-neighbour search algorithm for nonparametric classification; Steele and Patterson (2000) give formulae for exact calculation of bootstrap estimates of expected prediction error for nearest-neighbor classifiers; and Lin (2001) suggests a nonparametric classification rule for univariate data, based on the minimum Kolmogorov distance between two populations.

1.3. *Summary.* Section 2 presents our main results in the univariate, two-population case, where at least one of the densities is not close to zero. Section 3 suggests ways of removing the latter constraint; Section 4 treats empirical choice of bandwidth; Section 5 addresses generalizations to multiple and multivariate populations; and Section 6 outlines numerical properties. For the sake of brevity, most proofs are omitted, being available in a longer version of the paper, available online [Hall and Kang (2002)]. However, a brief account of the reasons for failure of leave-one-out methods is given in Section 7.

## 2. Classifying data from the body of a distribution: two-population case.

2.1. *Kernel-based classifiers.* Let the two populations have distributions $F$ and $G$, with respective densities $f$ and $g$. Let $0 < p < 1$ reflect the prior probability that a new, unclassified datum, $x$ say, lying in a given interval $\mathcal{I}$, is drawn from $F$. (To avoid degeneracy we assume throughout that $0 < p < 1$.) Denote by $\mathcal{A}_0$ the "ideal" algorithm that classifies $x$ as coming from $F$ or $G$ according as $\Delta(x) \equiv pf(x) - (1-p)g(x)$ is positive or negative, respectively. [We may make the classification arbitrarily if $\Delta(x)$ vanishes.] Among all measurable algorithms $\mathcal{A}$ for classification on $\mathcal{I}$, $\mathcal{A}_0$ is optimal in the sense of minimizing the Bayes risk

$$
\begin{aligned}
(2.1) \quad &\mathrm{err}_{\mathcal{A}}(f,g|\mathcal{I}) \\
&= p \int_{\mathcal{I}} P(x \text{ is classified by } \mathcal{A} \text{ as coming from } g) f(x)\,dx \\
&\quad + (1-p) \int_{\mathcal{I}} P(x \text{ is classified by } \mathcal{A} \text{ as coming from } f) g(x)\,dx.
\end{aligned}
$$

Optimality requires that prior probabilities for $F$ and $G$, restricted to $\mathcal{I}$, be precisely $p$ and $1-p$, respectively, although this assumption will not be a prerequisite for our main theoretical results.

Given training datasets $\mathcal{X} = \{X_1, \ldots, X_m\}$ and $\mathcal{Y} = \{Y_1, \ldots, Y_n\}$ drawn from $F$ and $G$, respectively, an empirical version of $\mathcal{A}_0$ may be based on nonparametric density estimators, $\hat{f}$ and $\hat{g}$ say, computed from $\mathcal{X}$ and $\mathcal{Y}$. Specifically, given a nonnegative kernel $K$ and bandwidths $h_1, h_2 > 0$, let

$$(2.2) \quad \hat{f}(x) = \frac{1}{mh_1} \sum_{i=1}^{m} K\left(\frac{x - X_i}{h_1}\right), \qquad \hat{g}(x) = \frac{1}{nh_2} \sum_{i=1}^{n} K\left(\frac{x - Y_i}{h_2}\right),$$

and let $\mathcal{A}_1$ be the rule that classifies $x$ as coming from $F$ or $G$, according as $\hat{\Delta}(x) \equiv p\hat{f}(x) - (1-p)\hat{g}(x)$ is positive or negative, respectively.

Classification can be made arbitrarily if $\hat{\Delta}(x) = 0$. However, in this case a distinction should be drawn between cases where at least one of $\hat{f}(x)$ and $\hat{g}(x)$ is nonzero and where $\hat{f}(x)$ and $\hat{g}(x)$ both vanish. In the latter setting classification can be more prone to error. An alternative algorithm, not employing arbitrary choice, will be suggested in Section 3.

2.2. *Main results.* We shall assume the following:

(2.3) $m/n$ is bounded away from zero and infinity as $n \to \infty$;

(2.4) $f$ and $g$ have two continuous derivatives and are bounded away from zero in an open interval containing $\mathcal{I}$;

(2.5) $\Delta$ vanishes at just $\nu \geq 1$ points, $y_1, \ldots, y_\nu$, in $\mathcal{I}$, all of them interior points and at each of which $\Delta'(y_j) \neq 0$;

(2.6) $K$ is a bounded, symmetric and compactly supported probability density;

(2.7) for $j = 1$ and $2$, $h_j = h_j(n) \asymp n^{-\rho}$ as $n \to \infty$, where $0 < \rho < 1$.

The notation $a(n) \asymp b(n)$ means that the ratio of left- and right-hand sides is bounded away from zero and infinity as $n \to \infty$. The equivalence of bandwidth sizes which (2.7) entails is not strictly necessary, but since optimal bandwidths satisfy (2.7), then it is imposed without loss of generality. Put $h = n^{-\rho}$, where $\rho$ is as in (2.7).

Our proof of Theorem 2.1, stated below, needs only two (or four, in the case of the second half of the theorem) continuous derivatives of $f$ and $g$ in neighborhoods of a cross-over point, together with continuity of $f$ and $g$ in an open interval $\mathcal{I}_{\mathrm{op}}$ containing $\mathcal{I}$, as asked by (2.4). However, (2.4) is a standard condition when analyzing performance of second-order density estimators, and two bounded derivatives are required for the minimax results of Marron (1983).

THEOREM 2.1. *Assume $0 < p < 1$ and $\mathcal{I}$ is a compact interval, and that* (2.3)–(2.7) *hold. Then,*

$$
(2.8) \quad \begin{aligned} &\mathrm{err}_{\mathcal{A}_1}(f,g|\mathcal{I}) - \mathrm{err}_{\mathcal{A}_0}(f,g|\mathcal{I}) \\ &= \tfrac{1}{2} \sum_{j=1}^{\nu} |\Delta'(y_j)|^{-1} E\{p\hat{f}(y_j) - (1-p)\hat{g}(y_j)\}^2 + o\{(nh)^{-1} + h^4\}. \end{aligned}
$$

*If in addition $\nu = 1$, $f''(y_1)g''(y_1) > 0$,*

$$
(2.9) \quad \frac{h_2}{h_1} = \left\{ \frac{pf''(y_1)}{(1-p)g''(y_1)} \right\}^{1/2} + o(h^2),
$$

*and $f$ and $g$ each have four continuous derivatives in a neighborhood of $y_1$, then* (2.8) *continues to hold if the remainder there is replaced by $o\{(nh)^{-1} + h^8\}$.*

Chanda and Ruymgaart (1989) give a version of (2.8) in cases where $g$ differs from $f$ only in location, and tails are controlled by specific decay assumptions. Result (2.8) is specific to kernel-based Bayes classifiers. Indeed, as we noted in Section 1.2, Cover (1968) has shown that much faster rates are possible for nearest-neighbor classifiers, for which the asymptotic risk usually dominates the Bayes risk $\mathrm{err}_{\mathcal{A}_0}(f,g|\mathcal{I})$.

An alternative algorithm is that suggested by Stoller (1954), and involves classifying a new data value $x$ as coming from $f$ if $x \leq \arg\max(m\widehat{F} - n\widehat{G})$, where $\widehat{F}$ and $\widehat{G}$ are the empirical distribution functions computed from $\mathcal{X}$ and $\mathcal{Y}$, respectively. Here the classification probability, for data in $\mathcal{I}$, converges to $\mathrm{err}_{\mathcal{A}_0}(f,g|\mathcal{I})$, but only at rate $O_p(n^{-1/2})$.

2.3. *Implications of Theorem* 2.1. The expansion at (2.8) may be refined to

$$(2.10) \quad \mathrm{err}_{\mathcal{A}_1}(f,g|\mathcal{I}) - \mathrm{err}_{\mathcal{A}_0}(f,g|\mathcal{I}) = B_1(nh)^{-1} + B_2 h^4 + o\{(nh)^{-1} + h^4\},$$

where $B_1$ and $B_2$ are both functions of $H_1 = h_1/h$ and $H_2 = h_2/h$, and, explicitly,

$$
(2.11) \quad \begin{aligned} B_1 &= \tfrac{1}{2}\kappa \sum_{j=1}^{\nu} |\Delta'(y_j)|^{-1} \{(rH_1)^{-1} p^2 f(y_j) + H_2^{-1}(1-p)^2 g(y_j)\}, \\ B_2 &= \tfrac{1}{8}\kappa_2^2 \sum_{j=1}^{\nu} |\Delta'(y_j)|^{-1} \{H_1^2 p f''(y_j) - H_2^2 (1-p) g''(y_j)\}^2, \end{aligned}
$$

with $\kappa = \int K^2$, $\kappa_j = \int u^j K(u)\, du$ and $r = m/n$. Result (2.10) implies that the optimal bandwidth is of size $n^{-1/5}$ (i.e., $\rho = 1/5$), and that optimal values of the constants $H_1$ and $H_2$ are obtained by minimizing $B_1 + B_2$, unless it should be possible to render $B_2 = 0$ by some positive, nonzero choice of $H_1$ and $H_2$.

If $\nu = 1$, then $B_2 = 0$ is possible (for positive $H_1$ and $H_2$) if and only if $f''(y_1)$ and $g''(y_1)$ are of the same sign; that is, the densities at the point $y_1$ where $pf$ and $(1-p)g$ cross are either both locally concave or both locally convex. Assuming this to be the case, and choosing $h_1$ and $h_2$ as at (2.9), we may show from (2.8) (with $h^8$ instead of $h^4$ in the remainder) that, instead of (2.10),

$$(2.12) \quad \text{err}_{\mathcal{A}_1}(f, g | \mathcal{I}) - \text{err}_{\mathcal{A}_0}(f, g | \mathcal{I}) = B_3(nh)^{-1} + B_4 h^8 + o\{(nh)^{-1} + h^8\},$$

where, defining $R = pf''(y_1)/(1-p)g''(y_1)$, we have

$$(2.13) \quad \begin{aligned} B_3 &= \frac{\kappa}{2H_1} |\Delta'(y_1)|^{-1} \{r^{-1} p^2 f(y_1) + R^{-1/2}(1-p)^2 g(y_1)\}, \\ B_4 &= \frac{\kappa_4^2 H_1^8}{1152} |\Delta'(y_1)|^{-1} \{pf^{(4)}(y_1) - R^2(1-p)g^{(4)}(y_1)\}^2. \end{aligned}$$

Result (2.12) implies that the optimal bandwidth is now of size $n^{-1/9}$ (i.e., $\rho = 1/9$), and that the optimal constant $H_1$ is obtained by minimizing $B_3 + B_4$.

There is, of course, a possibility that the factor $T(f, g) \equiv pf^{(4)}(y_1) - R^2(1-p)g^{(4)}(y_1)$ appearing in the definition of $B_4$ vanishes. In this case the term in $B_4 h^8$ at (2.12) should be replaced by one in $h^{12}$, and the remainder replaced by $o\{(nh)^{-1} + h^{12}\}$, provided $f$ and $g$ have continuous derivatives of order 6 in a neighborhood of $y_1$. However, since $T(f, g)$ is a particularly unusual functional of second and fourth derivatives of two distinct densities, then it is unlikely that in practice $T(f, g) = 0$.

In summary, excepting pathological cases that can be expected to arise only rarely, the optimal bandwidths for classification when $\nu \geq 2$ are $h_j^0 = H_j n^{-1/5}$, where $H_1, H_2 > 0$ are chosen to minimize

$$(2.14) \quad \begin{aligned} \sum_{j=1}^{\nu} |\Delta'(y_j)|^{-1} [&\kappa\{(rH_1)^{-1} p^2 f(y_j) + H_2^{-1}(1-p)^2 g(y_j)\} \\ &+ \tfrac{1}{4} \kappa_2^2 \{H_1^2 pf''(y_j) - H_2^2(1-p)g''(y_j)\}^2]. \end{aligned}$$

If $f''(y_1)g''(y_1) < 0$, then this prescription is also valid for $\nu = 1$. However, if $\nu = 1$ and $f''(y_1)g''(y_1) > 0$, then, excepting pathological cases where $T(f, g) = 0$, the optimal bandwidths are $h_1^0 = H_1 n^{-1/9}$ and $h_2^0 = H_2 n^{-1/9} = H_1 R^{1/2} n^{-1/9}$, where $H_1 > 0$ minimizes

$$(2.15) \quad \begin{aligned} &\frac{\kappa}{H_1} \{r^{-1} p^2 f(y_1) + R^{-1/2}(1-p)^2 g(y_1)\} \\ &+ \frac{\kappa_4^2 H_1^8}{576} \{pf^{(4)}(y_1) - R^2(1-p)g^{(4)}(y_1)\}^2. \end{aligned}$$

An extreme case is that where $\Delta$ is smooth and vanishes over a "plate," that is, a nondegenerate interval $\mathcal{J} = [a,b]$. Then, each derivative of $\Delta$ which exists must vanish on $\mathcal{J}$. Therefore, if no discontinuities of derivatives enter into the determination of properties of $\Delta$, the problem of estimating the endpoints of $\mathcal{J}$ is essentially parametric. Provided there are no other points where $\Delta$ vanishes, then it may be shown that under appropriate regularity conditions, an empirical rule can get within $O(n^{-1})$ of $\text{err}_{\mathcal{A}_0}(f,g|\mathcal{I})$.

The setting where Bayes risk equals zero is sometimes addressed in the context of machine learning [see, e.g., Ehrenfeucht, Haussler, Kearns and Valiant (1989)]. Excluding the uninteresting degenerate case in which $p(1-p) = 0$, and pathological cases where the support of $f$ starts exactly at a point where that of $g$ ends (or vice versa), this setting entails $\Delta$ vanishing on a plate, as discussed in the previous paragraph. Therefore, its main implications are those that we have discussed previously.

In many circumstances the discussion of classification given following Theorem 2.1 applies in a general, global sense, to an empirical algorithm $\widehat{\mathcal{A}}$ applied to any new datum $x \in \mathbb{R}$, rather than only to the algorithm $\mathcal{A}_1$ restricted to $\mathcal{I}$. Details will be given in the next section.

## 3. Classification in the tails.

3.1. *Kernel-based classifiers.* We shall assume that the supports of both $f$ and $g$ are intervals, that neither density vanishes in the interior of its support, and that a classification rule is sought in the upper tail. In this instance our algorithm will be based on the assumption that, sufficiently far to the right, the tail of $f$ exceeds that of $g$, or vice versa. Formally, we ask that either $f(x) > g(x)$ for all $x \in (x_0, x_{\text{supp}})$, or $g(x) > f(x)$ for all $x \in (x_0, x_{\text{supp}})$, where $x_0$ is strictly less than the right-hand end, $x_{\text{supp}}$, of the support of $f$ or $g$, respectively; and we seek a means of classifying new data $x > x_0$. Of course, $x_{\text{supp}}$ may be infinite.

If $x > x_0$ and $\hat{f}(x) = \hat{g}(x) = 0$, let $\hat{x}$ denote the infimum of values of $y \leq x$ such that $\hat{f}(z) = \hat{g}(z) = 0$ for all $z \in [y, x]$. Our algorithm, to which we refer below as $\mathcal{A}_R$, where the subscript indicates the right-hand tail, consists of classifying $x$ as coming from $f$ or $g$, according, as $\hat{f}(\hat{x}-) > 0$ or $\hat{g}(\hat{x}-) > 0$. [With probability 1, exactly one of $\hat{f}(\hat{x}-)$ and $\hat{g}(\hat{x}-)$ will be nonzero.]

3.2. *Main results.* Theorem 3.1 below shows that the suboptimality level discussed in Section 2, that is, $O(n^{-(1-\rho)})$ where $\rho = 1/5$ or $1/9$, is preserved if the upper tail weights of $f$ and $g$ are sufficiently different. Theorem 3.2 demonstrates by example that if the tail weights are too close, then the level of suboptimality can be of strictly larger order than $n^{-(1-\rho)}$.

Next we give regularity conditions for Theorem 3.1. Writing $F$ and $G$ for the distributions corresponding to densities $f$ and $g$, respectively, we ask that:

(3.1) $K$ is a bounded, symmetric, compactly supported and Hölder continuous probability density;

(3.2) for $j = 1$ and $2$, $h_j = h_j(n) \asymp n^{-\rho}$ as $n \to \infty$, where $0 < \rho < 1$;

(3.3) $f$ and $g$ are continuous, and strictly decreasing in their upper tails;

(3.4) for a constant $A_1 > 0$ and all sufficiently large $x$, $A_1 f(x) > f(x - x^{-1})$;

(3.5) each $\varepsilon > 0$ and all sufficiently large $x$, $\varepsilon f(x) \geq g(x - x^{-1})$;

(3.6) $A_2 > 0$, for $a > \frac{2-\rho}{1-\rho}$ and all sufficiently large $x$, $1 - G(x) \leq A_2 f(x)^a$;

(3.7) $x^{(2-\rho)/\rho}\{1 - G(x)\} \to 0$ as $x \to \infty$.

Assumption (3.1) is satisfied by compactly supported kernels commonly used in practice, and, in particular, by the Epanechnikov, biweight and triweight kernels; condition (3.2) is satisfied by the optimal bandwidths discussed in Section 2; (3.3) asks that the tails of $f$ and $g$ be smooth and eventually decreasing; (3.4) asks that the tails of $f$ not decrease too rapidly, and is satisfied by the majority of distributions that have infinite tails to the right; (3.5) asks that $f$ eventually dominate $g$; (3.6) asserts that this domination is sufficiently great; and (3.7) holds if the lighter-tailed distribution $G$ has finite moment of order $(2 - \rho)/\rho$.

THEOREM 3.1. *If (3.1)–(3.7) hold, then for some $x_0 > 0$,*

$$\begin{aligned}(3.8) \quad &P\{\text{for each } x > x_0, \text{ one of the following two properties} \\ &\text{holds: (a) } p\hat{f}(x) > (1-p)\hat{g}(x), \text{ or (b) } \hat{f}(x) = \hat{g}(x) = 0, \\ &\hat{g}(y) = 0 \text{ for all } y > x, \hat{f}(\hat{x}-) > 0 \text{ and } \hat{g}(\hat{x}-) = 0\} \\ &= 1 - o\{(nh)^{-1}\}\end{aligned}$$

*as $n \to \infty$.*

Next we investigate an instance where $f$ and $g$ both have Pareto-type tails, but the tail weights are sufficiently similar for the algorithm $\mathcal{A}_R$ to have difficulty distinguishing between them. Specifically, assume that

(3.9) $$f(x) \sim ax^{-\alpha} \quad \text{and} \quad g(x) \sim bx^{-\beta} \qquad \text{as } x \to \infty,$$
$$\text{where } a, b > 0 \text{ and } 1 < \alpha < \beta < \alpha + 1 < \infty.$$

Let $\mathcal{A}_2 = \mathcal{A}_1 \cup \mathcal{A}_R$ denote the algorithm constructed by using $\mathcal{A}_1$ to classify $x$ if not both of $\hat{f}(x)$ and $\hat{g}(x)$ vanish, and using $\mathcal{A}_R$ otherwise.

THEOREM 3.2. *If* (3.1), (3.2) *and* (3.9) *hold, then for all sufficiently large* $x_0$,

$$(3.10) \quad nh \int_{x_0}^{\infty} P(x \text{ is classified by } \mathcal{A}_2 \text{ as coming from } g) f(x) \, dx \to \infty$$

*as* $n \to \infty$.

3.3. *Implications of Theorems* 3.1 *and* 3.2. An immediate consequence of Theorem 3.1 is that if (3.1)–(3.7) hold, then the probability that, uniformly in new data $x$ on $[x_0, \infty)$, $\mathcal{A}_2$ is equivalent to classifying in the optimal way using $\mathcal{A}_0$, equals $1 - o\{(nh)^{-1}\}$. Therefore, taking the classification interval to be $\mathcal{I} = [x_0, \infty)$, we deduce that

$$(3.11) \quad \operatorname{err}_{\mathcal{A}_2}(f, g | \mathcal{I}) - \operatorname{err}_{\mathcal{A}_0}(f, g | \mathcal{I}) = o\{(nh)^{-1}\}$$

as $n \to \infty$. The left-hand side of (3.11) is of course nonnegative; it represents the Bayes risk for an empirical classification rule, minus the risk for the optimal rule.

There is, of course, a version of $\mathcal{A}_R$ for the left-hand tail; call it $\mathcal{A}_L$. Let $\widehat{\mathcal{A}}$ denote the algorithm that classifies $x$ using $\mathcal{A}_1$ if $\hat{f}(x)$ and $\hat{g}(x)$ do not both vanish, or using $\mathcal{A}_R$ if $\hat{f}(x) = \hat{g}(x) = 0$ and $x$ lies to the right of the median of $\mathcal{X} \cup \mathcal{Y}$, or using $\mathcal{A}_L$ otherwise. (Our choice of the median is arbitrary.) Assume $f$ and $g$ are continuous on the real line, that the supports of $f$ and $g$ are intervals, that neither density vanishes at any point in the interior of its support, that the conditions of Theorem 2.1 hold on any compact interval $\mathcal{I}$ that is interior to the intersection of the supports, that the conditions of Theorem 3.1 (possibly with $f$ and $g$ interchanged) hold to the right, and that the analogous conditions hold to the left. Then in either tail, either $f$ or $g$ dominates the other, and so there can be only a finite number of points ($\nu$, say) at which the graphs of $pf$ and $(1-p)g$ cross.

In these circumstances we may deduce from Theorems 2.1 and 3.1 that the expansions of classification error described in Theorem 2.1 hold for the algorithm $\widehat{\mathcal{A}}$ applied to classification on the whole real line $\mathbb{R}$:

$$(3.12) \quad \begin{aligned} & \operatorname{err}_{\widehat{\mathcal{A}}}(f, g | \mathbb{R}) - \operatorname{err}_{\mathcal{A}_0}(f, g | \mathbb{R}) \\ & = \tfrac{1}{2} \sum_{j=1}^{\nu} |\Delta'(y_j)|^{-1} E\{p\hat{f}(y_j) - (1-p)\hat{g}(y_j)\}^2 + o\{(nh)^{-1} + h^4\}. \end{aligned}$$

The remainder term here can be sharpened to $o\{(nh)^{-1} + h^8\}$ if the conditions of the second part of Theorem 2.1 apply, in particular, if $h_1$ and $h_2$ satisfy (2.9).

In view of these results, the discussion of optimality given following Theorem 2.1 applies to the present general, global setting, where $\widehat{\mathcal{A}}$ is used to classify any real-valued datum $x$. The asymptotically optimal bandwidths are either $h_j^0 = H_j n^{-1/5}$ or $h_j^0 = H_j n^{-1/9}$, where $(H_1, H_2)$ minimizes either (2.14) or (2.15), respectively, and $H_2 = R^{1/2} H_1$ in the latter case.

It may be proved from (3.9) that if $x_0$ is sufficiently large, then (3.11) fails. Therefore, if the bandwidths $h_1$ and $h_2$ are chosen so as to minimize the inherent additional classification error in the body of the distribution, relative to the optimal algorithm $\mathcal{A}_0$, this performance will not be reflected when using $\mathcal{A}_2$ to classify data in the tails. If (3.9) holds, then the additional error introduced by the difficulty of classifying data in the tails is so large as to dominate the relatively low levels of error (in comparison with $\mathcal{A}_0$) experienced elsewhere.

The rate of divergence in (3.10) can be arbitrarily slow, in the sense that for any given $\varepsilon > 0$ there exist densities $f$ and $g$ satisfying (3.9) and for which the left-hand side of (3.10) diverges to infinity more slowly than $n^\varepsilon$, as $n \to \infty$.

Work of Chanda and Ruymgaart (1989) provides some further detail related to Theorem 3.2. Addressing the case where $f$ and $g$ differ only in location, and the density tails decrease like $x^{-\gamma}$ as $x \to \infty$, Chanda and Ruymgaart show that the difference between the error of the empirical classifier and its asymptotic limit is of size $(nh)^{-\gamma/(\gamma+2)}$. Moreover, if the density tails decrease like $e^{-x^\gamma}$, then the rate $O(n^{-4/5})$ is possible if $\gamma > 1$, although a slower rate occurs if $\gamma \leq 1$.

## 4. Empirical choice of bandwidth.

4.1. *Discussion of methods.* We could compute bandwidths by constructing empirical approximations to the functions appearing in (2.14) and (2.15), finding the minima of empirical forms of those expressions and substituting the resulting values into formulae for theoretically optimal bandwidths. However, this technique is awkward to use, since it requires explicitly working out how many times the graphs of $pf$ and $(1-p)g$ cross and where the crossings take place. This calls for technology similar to bump hunting methods. The relative complexity of that approach motivates alternative, more implicit techniques for bandwidth selection. One possibility is cross-validation, which at first sight seems very attractive.

A cross-validation method for choosing bandwidth is as follows. Let $\hat{f}_{-i}$ and $\hat{g}_{-i}$ denote the respective versions of $\hat{f}$ and $\hat{g}$, defined at (2.2), that are obtained through computing the latter estimators from the leave-one-out datasets $\mathcal{X}_i = \mathcal{X} \setminus \{X_i\}$ and $\mathcal{Y}_i = \mathcal{Y} \setminus \{Y_i\}$, respectively. (We continue to use respective bandwidths $h_1$ and $h_2$.) Put $\widehat{\Delta}_{f,-i} = p\hat{f}_{-i} - (1-p)\hat{g}$, $\widehat{\Delta}_{g,-i} =$

$p\hat{f} - (1-p)\hat{g}_{-i}$ and

$$\widetilde{\mathrm{err}}_{\mathcal{A}_1}(h_1, h_2) = \frac{p}{m} \sum_{i=1}^{m} I\{\widehat{\Delta}_{f,-i}(X_i) < 0, X_i \in \mathcal{I}\}$$
$$+ \frac{1-p}{n} \sum_{i=1}^{n} I\{\widehat{\Delta}_{g,-i}(Y_i) > 0, Y_i \in \mathcal{I}\}.$$
(4.1)

One might choose $(h_1, h_2) = (\hat{h}_1, \hat{h}_2)$ to minimize $\widetilde{\mathrm{err}}_{\mathcal{A}_1}(h_1, h_2)$. The latter may be viewed as an empirical approximation to $\mathrm{err}_{\mathcal{A}_1}(f, g|\mathcal{I})$. However, this approach performs poorly in both theory and practice, and, in particular, does not accurately estimate, in the sense of relative consistency, the value of $(h_1, h_2)$ that minimizes $\mathrm{err}_{\mathcal{A}_1}(f, g|\mathcal{I})$. See Section 7 for details.

A second, more effective approach, which we shall consider in detail, is based on using the bootstrap to estimate $\mathrm{err}_{\mathcal{A}_1}(f, g|\mathcal{I})$ and, thereby, to select the optimal bandwidths. Specifically, let $\tilde{f}$ and $\tilde{g}$ be the versions of $\hat{f}$ and $\hat{g}$, defined at (2.2), that arise if we use respective bandwidths $h_3$ and $h_4$ (instead of $h_1$ and $h_2$). Conditional on $\mathcal{X}$ (or on $\mathcal{Y}$), draw $m$ data $\mathcal{X}^* = \{X_1^*, \ldots, X_m^*\}$ independently and uniformly from the distribution with density $\tilde{f}$ (or, resp., $n$ data $\mathcal{Y}^* = \{Y_1^*, \ldots, Y_n^*\}$ independently and uniformly from the distribution with density $\tilde{g}$), and let

$$\hat{f}^*(x) = \frac{1}{mh_1} \sum_{j=1}^{m} K\left(\frac{x - X_j^*}{h_1}\right), \qquad \hat{g}^*(x) = \frac{1}{nh_2} \sum_{j=1}^{n} K\left(\frac{x - Y_j^*}{h_2}\right).$$

Put $\widehat{\Delta}^*(x) = p\hat{f}^*(x) - (1-p)\hat{g}^*(x)$ and

$$\widehat{\mathrm{err}}_{\mathcal{A}_1}(h_1, h_2) = p \int P\{\widehat{\Delta}^*(x) < 0 | \mathcal{X} \cup \mathcal{Y}\} \tilde{f}(x)\, dx$$
$$+ (1-p) \int P\{\widehat{\Delta}^*(x) > 0 | \mathcal{X} \cup \mathcal{Y}\} \tilde{g}(x)\, dx.$$

Choose $(h_1, h_2) = (\hat{h}_1, \hat{h}_2)$ to minimize $\widehat{\mathrm{err}}_{\mathcal{A}_1}(h_1, h_2)$.

In the two respective cases we need to choose $h_3$ and $h_4$ so that the "pilot" density estimators $\tilde{f}$ and $\tilde{g}$ are able to consistently estimate second, or fourth, derivatives of $f$ and $g$. It is known from more conventional applications of curve estimation that this requires $h_3$ and $h_4$ to be of strictly larger order than $n^{-1/5}$ or $n^{-1/9}$, respectively. Therefore, we should choose $h_3$ and $h_4$ to both be of size $n^{-\sigma}$, where in the first regime $0 < \sigma < \frac{1}{5}$ and in the second $0 < \sigma < \frac{1}{9}$. Since taking $0 < \sigma < \frac{1}{9}$ covers both cases, then, for simplicity, we shall make that assumption in our theoretical results below. For the same reason we shall assume four derivatives of $f$ and $g$ in the neighborhood of each cross-over point, although in the case of the first regime only two derivatives are required.

4.2. *Main results.* We shall assume the following:

(4.2) $f$ and $g$ are continuously differentiable and are bounded away from zero on an open interval containing $\mathcal{I}$; $\Delta$ vanishes at just $\nu$ points, $y_1, \ldots, y_\nu$, in $\mathcal{I}$, all of them interior points and at each of which $\Delta'(y_j) \neq 0$; and $f$ and $g$ each have four continuous derivatives in neighborhoods of each $y_j$;

(4.3) either $\nu \geq 1$, $f''(y_j)g''(y_j) \neq 0$ for at least one $j$, and the $\nu$ equations $pf''(y_j) - R(1-p)g''(y_j) = 0$ do not have a simultaneous solution $R > 0$; or $\nu = 1$ and $f''(y_j)g''(y_j) > 0$, in which case a solution exists;

(4.4) $n/n$ is bounded away from zero and infinity as $n \to \infty$;

(4.5) $K$ is a compactly supported function with four Hölder continuous derivatives on the real line and satisfying $\int K = 1$;

(4.6) for $j = 3$ and 4, $h_j = h_j(n) \asymp n^{-\sigma}$ as $n \to \infty$, where $0 < \sigma < 1/15$.

Condition (4.3) implies that one or other of the two main regimes of behavior of $h_1^0$ and $h_2^0$ obtains. If $\rho = \frac{1}{5}$ or $\frac{1}{9}$ in the two respective cases, then the optimal bandwidths are $h_j^0 \sim H_j n^{-\rho}$ for $j = 1, 2$, where $H_1$ and $H_2$ are positive constants.

Given $0 < c_1 < \frac{1}{9} < \frac{1}{5} < c_2 < 1$, let $(h_1, h_2) = (\hat{h}_1, \hat{h}_2)$ denote the bandwidth pair that minimizes $\widehat{\text{err}}_{\mathcal{A}_1}(h_1, h_2)$ over $(h_1, h_2)$ such that $n^{-c_2} \leq h_1, h_2 \leq n^{-c_1}$. The theorem below shows that each empirical bandwidth $\hat{h}_j$ is asymptotic to its asymptotically optimal counterpart $h_j^0$. In addition, if a sufficiently high-order kernel is used to estimate $\tilde{f}$ and $\tilde{g}$, then an empirical form of (2.9) holds.

THEOREM 4.1. *Assume $0 < p < 1$ and $\mathcal{I}$ is a compact interval, and that (4.2)–(4.6) hold. Then, for $j = 1$ and 2, $\hat{h}_j/h_j^0 \to 1$ in probability as $n \to \infty$. Furthermore, if $K$ is of order $r$, meaning that $\int u^j K(u)\,du = 0$ for $1 \leq j \leq r-1$, if $r > 2/(5\sigma)$, if the second part of (4.3) obtains, and if $f$ and $g$ have $r + 2$ bounded derivatives in a neighborhood of $y_1$, then the following empirical form of (2.9) holds:*

$$(4.7) \qquad \frac{\hat{h}_2}{\hat{h}_1} = \left\{ \frac{pf''(y_1)}{(1-p)g''(y_1)} \right\}^{1/2} + o_p(n^{-2\rho}).$$

## 5. Multiple or multivariate populations.

5.1. *Multiple univariate populations.* Suppose there are $N$ distributions, $F_1, \ldots, F_N$ say, with respective densities $f_1, \ldots, f_N$ and prior probabilities $p_1, \ldots, p_N$, where $\sum_j p_j = 1$. Let $\mathcal{A}$ denote a general algorithm for classifying data in a given interval $\mathcal{I}$. The "ideal" algorithm which minimizes the Bayes risk

$$\mathrm{err}_{\mathcal{A}}(f_1,\ldots,f_N|\mathcal{I})$$
$$= \sum_{j=1}^{N} p_j \int_{\mathcal{I}} P(x \text{ is not classified by } \mathcal{A} \text{ as coming from } f_j) f_j(x)\,dx,$$

is the classification rule $\mathcal{A}_0$ which declares $x$ to have come from $f_j$ if $p_j f_j(x) = \max_k \{p_k f_k(x)\}$. (Ties may be broken at random.) Here it is assumed that the prior probabilities for $f_1,\ldots,f_N$, restricted to $\mathcal{I}$, are $p_1,\ldots,p_N$, respectively.

Assume that for each $1 \leq j \leq N$, we have access to a sample $X_{j1},\ldots,X_{jn_j}$ of independent and identically distributed data drawn from distribution $F_j$. Assume the samples are themselves independent. Construct the density estimator

$$\hat{f}_j(x) = \frac{1}{n_j h_j} \sum_{i=1}^{n_j} K\left(\frac{x - X_{ji}}{h_j}\right),$$

where $h_j$ is a bandwidth. Let $\mathcal{A}_1$ denote the empirical algorithm which declares $x$ to have come from $f_j$ if and only if $p_j \hat{f}_j = \max_k \{p_k \hat{f}_k(x)\}$. (Breaking ties at random in this rule has no effect on our asymptotic results, provided $\max_j p_j f_j$ is bounded away from zero on $\mathcal{I}$.)

Let $\mathcal{I}$ denote a compact interval, and assume $\max_j p_j f_j$ is bounded away from zero in an open interval containing $\mathcal{I}$; that $\Delta_{ij} \equiv p_i f_i - p_j f_j$ vanishes only at discrete interior points $y_{ijk}$ of $\mathcal{I}$, where $1 \leq k \leq \nu_{ij}$ and $\Delta'_{ij}(y_{ijk}) \neq 0$; that these points are distinct, in the sense that $y_{i_1 j_1 k_1} = y_{i_2 j_2 k_2}$ implies $\{i_1, j_1\} = \{i_2, j_2\}$ and $k_1 = k_2$; that $n_1 \to \infty$ and each ratio $n_j/n_k$ is bounded; that for each $1 \leq j \leq N$, $h_j = h_j(n_1) \asymp n_1^{-1/5}$ as $n_1 \to \infty$; and that other conditions, for example, on the smoothness of each $f_j$, are analogous to those in Section 2. Put $n = n_1$, $H_j = n^{1/5} h_j$ and $r_j = n_j/n$, let $\kappa$ and $\kappa_2$ be as in Section 2, and define

$$
\begin{aligned}
T(H_1,\ldots,H_N) &\\
= \frac{\kappa}{4} \sum_{i \neq j} \sum_{k=1}^{\nu_{ij}} |\Delta'(y_{ijk})|^{-1} &\{(r_i H_i)^{-1} p_i^2 f_i(y_{ijk}) \\
& + (r_j H_j)^{-1} p_j^2 f_j(y_{ijk})\} \\
+ \frac{\kappa_2^2}{16} \sum_{i \neq j} \sum_{k=1}^{\nu_{ij}} |\Delta'(y_{ijk})|^{-1} &\{H_i^2 p_i f_i''(y_{ijk}) - H_j^2 p_j f_j''(y_{ijk})\}^2.
\end{aligned}
$$
(5.1)

Then the following analogues of (2.8) and (2.10) may be derived:

$$
\begin{aligned}
&\mathrm{err}_{\mathcal{A}_1}(f_1,\ldots,f_N|\mathcal{I}) - \mathrm{err}_{\mathcal{A}_0}(f_1,\ldots,f_N|\mathcal{I}) \\
&= \tfrac{1}{4}\sum_{i\neq j}\sum \sum_{k=1}^{\nu_{ij}} |\Delta'(y_{ijk})|^{-1} E\{p_i \hat{f}_i(y_{ijk}) - p_j \hat{f}_j(y_{ijk})\}^2 + o(n^{-4/5}) \\
&= T(H_1,\ldots,H_N) n^{-4/5} + o(n^{-4/5}).
\end{aligned}
\tag{5.2}
$$

5.2. *Implications of* (5.2). Our assumptions imply that no three graphs of the functions $p_i f_i$ cross at a single point $y \in \mathcal{I}$, and, indeed, (5.2) fails in such cases. Although those cases might be considered rare, Fukunaga and Flick (1984) show that they can arise.

The context directly addressed by (5.2) is that where it is impossible to choose $H_1,\ldots,H_N > 0$ such that $(H_i/H_j)^2 = p_j f_j''(y_{ijk})/\{p_i f_i''(y_{ijk})\}$ for each triple of indices $(i,j,k)$ such that $p_i f_i$ and $p_j f_j$ cross at some point $y_{ijk} \in \mathcal{I}$. For example, this can be because $f_j''(y_{ijk})f_i''(y_{ijk}) < 0$ for some $(i,j,k)$, or because for some pair $(i,j)$ the ratio $f_i''(y_{ijk})/f_j''(y_{ijk})$ varies with $k$. Here the optimal rate of convergence to zero of the difference in Bayes risk is $n^{-4/5}$, and its minimal size is obtained by choosing $h_j = H_j n^{-1/5}$, where $H_1,\ldots,H_N$ minimizes $T(H_1,\ldots,H_N)$ at (5.1).

Consider next the case where there is only one nonzero value of $\nu_{ij}$, and it equals 1. Here the optimal algorithm $\mathcal{A}_0$ reduces to distinguishing between just two densities, $f_i$ and $f_j$ say. The empirical algorithm $\mathcal{A}_1$ also effectively reduces to a two-population one, where the convergence rate can be either $n^{-4/5}$ or $n^{-8/9}$. Since this case has already been discussed in Section 2, then there is no need to treat it further.

There are, however, nonpathological instances where $\sum_{i<j} \nu_{ij} > 1$ and the convergence rate $n^{-8/9}$, rather than $n^{-4/5}$, obtains. Consider, for example, the case where, for $1 \leq j \leq M$ (and $M < N$), the graph of $p_j f_j$ crosses that of $p_{j+1} f_{j+1}$ at a single point, $y_j$ say; and no other crossings of graphs occur within $\mathcal{I}$. If, at each crossing, the graphs are all locally concave or all locally convex, then, by choosing $H_1,\ldots,H_{M+1}$ such that $(H_j/H_{j+1})^2 = p_{j+1} f_{j+1}''(y_j)/\{p_j f_j''(y_j)\}$ for $1 \leq j \leq M$, we ensure that the bias contribution to $T(H_1,\ldots,H_N)$, that is, the second term in (5.1), vanishes identically. In this case the faster convergence rate of $n^{-8/9}$ can be obtained by choosing $h = H_j n^{-1/9}$ throughout. (Choice of $h_j$ for $j > M+2$ is relatively unimportant, since the corresponding densities do not cross any other density in $\mathcal{I}$. Nevertheless, taking $h_j = n^{-1/9}$ is adequate.) There are many related examples of this type.

5.3. *Multivariate populations.* Let $f$ and $g$ be densities of $d$-variate distributions $F$ and $G$, respectively, where $d \geq 1$. We assume classification is conducted for new data $x$ coming from a region $\mathcal{R}$, which here plays the role of the interval $\mathcal{I}$ in Section 2. The empirical rule $\mathcal{A}_1$ classifies $x$ as coming from $F$ or $G$, according, as $p\hat{f}(x) - (1-p)\hat{g}(x)$ is positive or negative, where on the present occasion,

$$\hat{f}(x) = \frac{1}{mh_1^d} \sum_{i=1}^{m} K\left(\frac{x - X_i}{h_1}\right), \qquad \hat{g}(x) = \frac{1}{nh_2^d} \sum_{i=1}^{n} K\left(\frac{x - Y_i}{h_2}\right),$$

$\mathcal{X} = \{X_1, \ldots, X_m\}$ and $\mathcal{Y} = \{Y_1, \ldots, Y_n\}$ are training datasets drawn from $F$ and $G$, respectively, $h_1$ and $h_2$ are bandwidths, and $K$ is a bounded, spherically symmetric and compactly supported probability density.

The classification rule $\mathcal{A}_0$ that minimizes Bayes risk amounts to classifying $x$ as coming from $F$ or $G$ according as $\Delta(x) > 0$ or $< 0$, where $\Delta = pf - (1-p)g$. Let $\mathcal{C}$ denote that part of the set $\{y : \Delta(y) = 0\}$ which lies in $\mathcal{R}$, and write $\theta(y)$ for the vector of first derivatives of $\Delta$ at $y$. In place of (2.4) and (2.5), we assume that $f$ and $g$ have two continuous derivatives, and are bounded away from zero, in an open set containing $\mathcal{R}$, and the function $\theta$ does not vanish on $\mathcal{C}$. Take each $h_j$ to be of size $n^{-1/(d+4)}$. Then it may be proved that

$$\begin{aligned}(5.3) \quad &\operatorname{err}_{\mathcal{A}_1}(f, g | \mathcal{R}) - \operatorname{err}_{\mathcal{A}_0}(f, g | \mathcal{R}) \\ &= \tfrac{1}{2} \int_{\mathcal{C}} \|\theta(y)\|^{-1} E\{p\hat{f}(y) - (1-p)\hat{g}(y)\}^2 \, dy + o(n^{-4/(d+4)}).\end{aligned}$$

Holmström and Klemelä (1992) report the results of numerical experiments on kernel-based classification in the multivariate case. They provide no theory, however.

5.4. *Implications of* (5.3). Taking $h_j = H_j n^{-1/5}$, Taylor expansion of the right-hand side of (5.3) may be shown to give

$$\operatorname{err}_{\mathcal{A}_1}(f, g | \mathcal{R}) - \operatorname{err}_{\mathcal{A}_0}(f, g | \mathcal{R}) = B(H_1, H_2) n^{-4/(d+4)} + o(n^{-4/(d+4)}),$$

where the constant $B(H_1, H_2)$ vanishes for either finite or infinite $(H_1, H_2)$ only if $\nabla^2 f / \nabla^2 g$ is constant throughout $\mathcal{C}$, with $\nabla^2 \psi$ denoting the Laplacian. Therefore, in virtually all cases there exists an optimal pair $(H_1, H_2) = (H_1^0, H_2^0)$ which minimizes $B(H_1, H_2)$. Then the optimal bandwidths $h_j^0 = H_j^0 n^{-1/(d+4)}$ are of size $n^{-1/(d+4)}$, which is the same size that leads to minimization of mean squared error of $\hat{f}$ and $\hat{g}$ as pointwise estimators of $f$ and $g$.

**6. Numerical properties.** We summarize a simulation study addressing properties of the empirical bandwidth selector introduced in Section 4. Recall from Section 2 that there are two main classes of problems, respectively characterized by the property that the densities $f$ and $g$ intersect at a point where the curvatures have different signs or the same sign. Call these classes 1 and 2; they correspond to the optimal bandwidth being of size $n^{-1/5}$ or $n^{-1/9}$, respectively. We shall report results for two examples in each class. Throughout, the distribution with density $f$ was standard normal, $p = \frac{1}{2}$ and $m = n$. In the tails of the distributions, in any cases of ambiguity we classified using the method suggested in Section 3.

Classification was done on the entire real line, rather than on a compact interval as suggested in our theory. In the first examples, in each of classes 1 and 2 the densities cross one another at one point in the tails, in addition to a crossing in the "middle" of the distribution. However, the tail crossing point is so far out that, for the sample sizes we used, it has negligible impact on numerical results, and so the effective value of $\nu$ is 1. The actual value of $\nu$ is 1 for the second example in class 1. For the second example in class 2, $\nu = 2$. However, there is strong symmetry in this case, with the result that theoretical properties are essentially the same as they would be if $\nu$ were 1. Nevertheless, the existence of two crossing points creates potential hazards for our empirical bandwidth selector, which is why we treated this example.

In the first example in class 1, $g$ is the $N(-1.2, 0.6^2)$ density, the crossover occurs at $y_1 = -0.515$, and the curvatures there are $f''(y_1) = -0.255$ and $g''(y_1) = 0.281$. In the second example in class 1, $g$ is the density for the normal mixture

$$\tfrac{1}{5}N(0,1) + \tfrac{1}{5}N(1, (\tfrac{2}{3})^2) + \tfrac{3}{5}N(\tfrac{19}{12}, (\tfrac{5}{9})^2),$$

$y_1 = 0.707$, $f''(y_1) = -0.156$ and $g''(y_1) = 0.327$. In the first example in class 2, $g$ is the normal $N(1, 1)$ density, $y_1 = 0.5$ and $f''(y_1) = g''(y_1) = -0.264$. In the second example in class 2, $g$ is the Cauchy density $g(x) = \{\pi(1+x^2)\}^{-1}$, there are two crossover points $y_i = \pm 1.851$, and $f''(y_i) = 0.175$ and $g''(y_i) = 0.068$. Figure 1 illustrates the densities.

To implement the bootstrap method suggested in Section 4, we used the triweight kernel, $K(x) = (35/32)(1-x^2)^3$ for $|x| \leq 1$, and noted that the asymptotically optimal bandwidth for estimating $f^{(r)}$, in terms of minimizing mean integrated squared error, is

$$h = \left\{ \frac{(2r+1)R(K^{(r)})}{n\mu_2(K)^2 \int (f^{(r+2)})^2} \right\}^{1/(2r+5)},$$

where $R(L) = \int L^2$ and $\mu_2(L) = \int u^2 L(u)\,du$. When constructing estimators $\tilde{f}$ and $\tilde{g}$ mentioned in Section 4, we took $r = 4$ and chose $h_3$, $h_4$ using the above formula, but (employing a device that might be implemented in practice) replaced $f$ by the normal density with zero mean and variance estimated from the training data. In the case of the Cauchy density, however, estimating scale in this way is inappropriate, and so instead the normalized interquartile range was used:

$$\hat{\sigma}_{\text{IQR}} = \frac{\text{sample interquartile range}}{\Phi^{-1}(0.75) - \Phi^{-1}(0.25)},$$

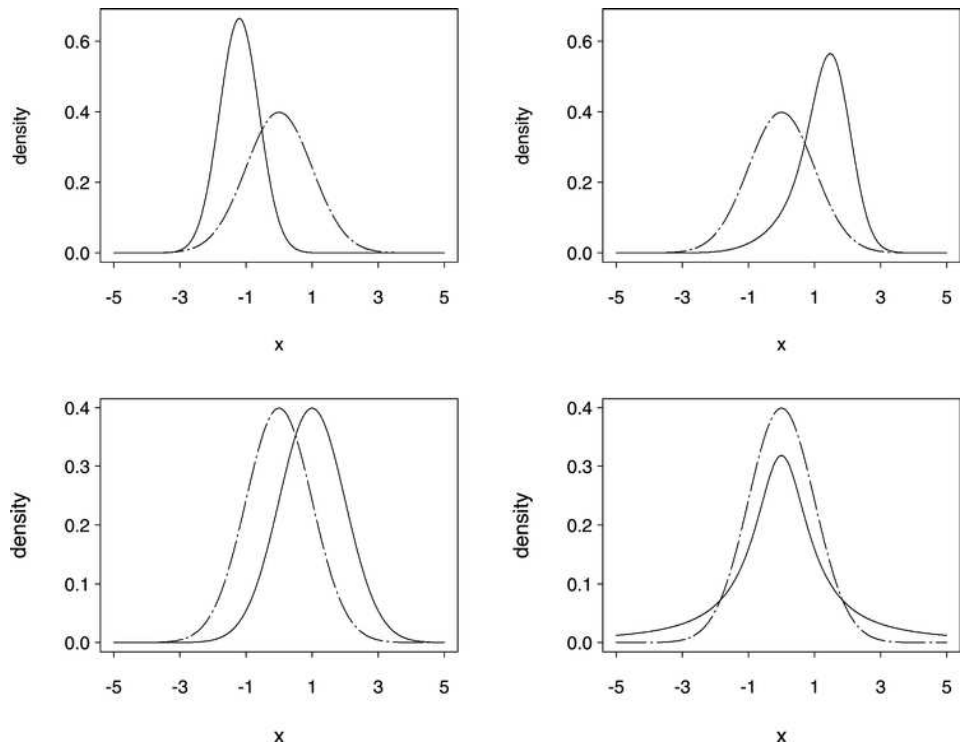

Fig. 1. *Densities used in simulation study. In each case the density $f(x) = (2\pi)^{-1/2} \exp(-\frac{1}{2}x^2)$ is indicated by the dot-dashed line, and the density $g$ by the unbroken line. The densities depicted in the two panels in the first and second rows correspond to those in the two examples in classes* 1 *and* 2, *respectively.*

where $\Phi^{-1}$ denotes the standard normal quartile function.

The probability $P\{\widehat{\Delta}^*(x) < 0 | \mathcal{X} \cup \mathcal{Y}\}$ needed to estimate $\widehat{\text{err}}_{\mathcal{A}_1}(h_1, h_2)$ was approximated using 100 bootstrap iterations. Minimization of $\widehat{\text{err}}_{\mathcal{A}_1}(h_1, h_2)$ over $(h_1, h_2)$ was conducted on a fine grid of bandwidths. We simulated 100 samples for each of 10 logarithmically equally spaced sample sizes from 20 to 200.

Let $(\hat{h}_1, \hat{h}_2)$ denote the empirical bandwidths obtained in this way. For each of the four distributions, and for $j = 1, 2$, we plotted $-\log \hat{h}_j$ against $\log n$. The results are given in Figures 2 and 3, which correspond to class 1 and class 2, respectively. In each figure, the two rows of panels give plots that correspond to the first and second density pairs, respectively, in that class; and the first and second columns of panels show (as black dots) the average values (over the 100 independent samples) of the points $(-\log \hat{h}_1, \log n)$ in the case of the left-hand panel, or $(-\log \hat{h}_2, \log n)$ for the right-hand panel. In each of the four panels in each figure, the unbroken line is the conventional least-squares regression line through these points. The dotted and dashed lines have slopes $\frac{1}{5}$ and $\frac{1}{9}$, respectively, with intercepts chosen so that each of these lines passes through the center of the least-squares regression line.

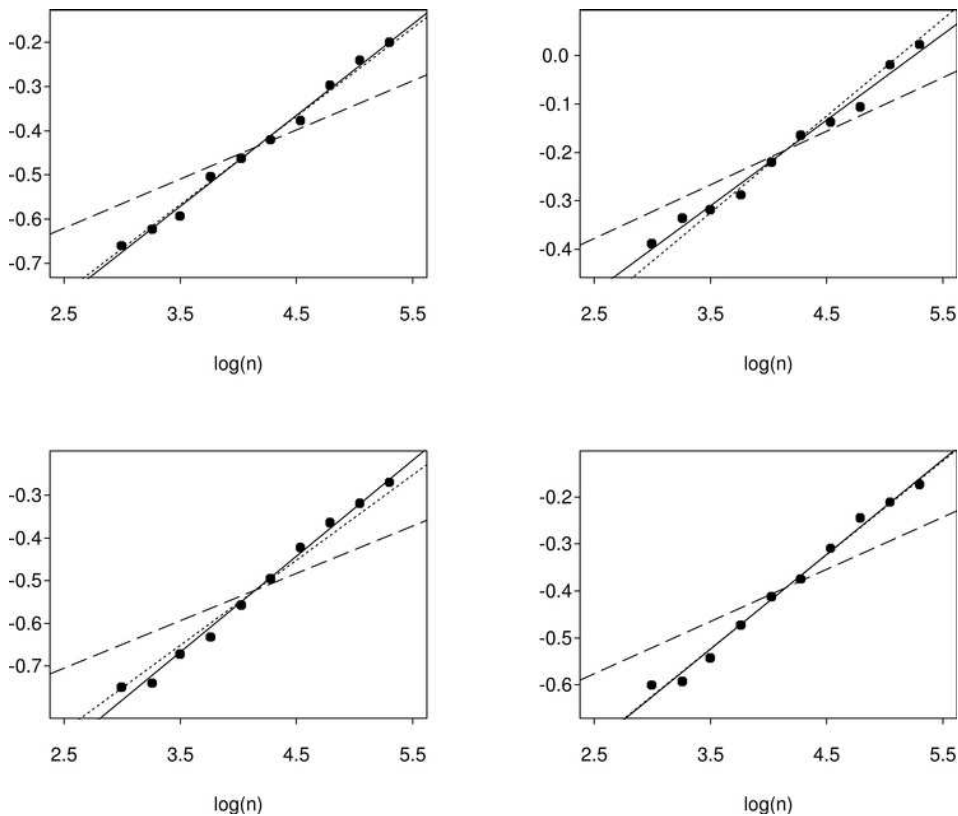

FIG. 2. *Plots for two examples in class* 1. *The two rows of panels show, respectively, simulation results for the two pairs of densities in class* 1, *that is, for the density pairs shown respectively in the first and second panels (in the first row) of Figure* 1. *In the jth column of each row the black dots show average values of* $(-\log \hat{h}_j, \log n)$, *computed as described in Section* 6. *The unbroken line is the conventional least-squares regression line through these points, and the dotted and dashed lines are drawn so that they have respective slopes* $\frac{1}{5}$ *and* $\frac{1}{9}$, *and pass through the center of the least-squares regression line.*

The main point to note from the figures is that in the case of density pairs from class 1, the slope of the least-squares regression line is very close to $\frac{1}{5}$ (see Figure 2), while for class 2 it is close to $\frac{1}{9}$ (see Figure 3). This, of course, reflects the theoretical results presented in Sections 3 and 4, where we showed that these particular slopes determine the optimal orders of bandwidth in the respective classes. The agreement between theory and numerical simulation is somewhat better in the case of class 1, but note that in the second class the numerical results clearly reflect the theory even in the Cauchy case.

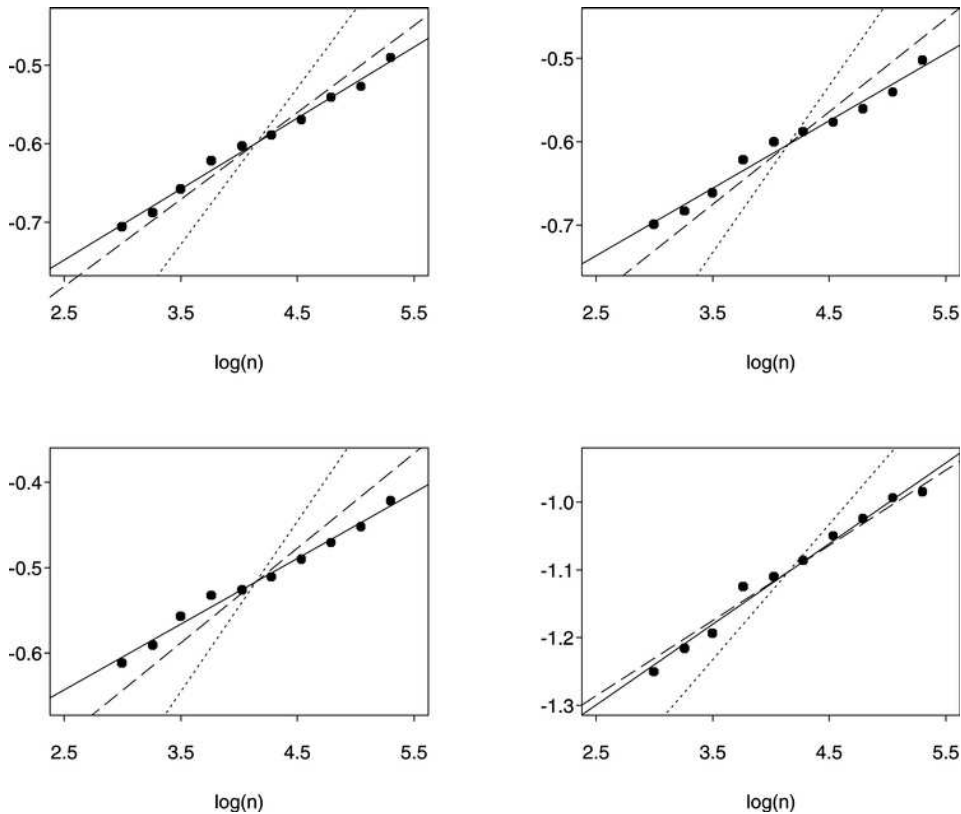

Fig. 3. *Plots for two examples in class* 2. *Details are as for Figure* 2, *except that the two rows of panels show results for the two pairs of densities in class* 2. *These density pairs are depicted in the first and second panels, respectively, in the last row of Figure* 1.

**7. Reasons for failure of $\widetilde{\mathrm{err}}_{\mathcal{A}_1}(h_1, h_2)$, at (4.1), to provide effective minimization of Bayes risk.** Failure occurs because the optimal bandwidths, discussed in Section 2.3, are determined by properties of mean squared error at isolated points, that is, the points where the graphs of $pf$ and $(1-p)g$ cross. See (2.8). Cross-validation does not accurately estimate mean squared error at a point, unless one averages over neighboring points in a sufficiently wide interval. See, for example, the modifications of cross-validation that are necessary when it is used for local, as distinct from global, bandwidth choice [Hall and Schucany (1989) and Mielniczuk, Sarda and Vieu (1989)]. The same sort of averaging is required here, too, and so the use of subsidiary smoothing parameters is necessary to overcome the failure of cross-validation. That substantially reduces the attractiveness of the method.

To appreciate these difficulties from a theoretical viewpoint, note that in order for the criterion defined at (4.1) to perform its function, it must equal $\mathrm{err}_{\mathcal{A}_1}(h_1, h_2)$, plus terms which either do not depend on $(h_1, h_2)$ or which depend on that quantity but are of smaller order than $\eta \equiv (mh_1)^{-1} + (nh_2)^{-1} + h_1^4 + h_2^4$. (We shall say that such terms are of "type T.") It is not difficult to see that this must be true of both series on the right-hand side of (4.1); there cannot, in general, be judicious cancellation between the two quantities. In particular,

$$S(h_1, h_2) \equiv \frac{1}{m} \sum_{i=1}^{m} I\{\widehat{\Delta}_{f,-i}(X_i) < 0, X_i \in \mathcal{I}\}$$

must equal $s(h_1, h_2) \equiv \int_{\mathcal{I}} P\{p\hat{f} < (1-p)\hat{g}\}f$, plus terms of type T; call this property $P_1$. We shall outline a theoretical argument showing that, in general, $P_1$ fails to hold.

For simplicity, let us take $p = \frac{1}{2}$, and $h_1$ and $h_2$ both to lie within the interval $\mathcal{H} = [n^{-1/5}C_1, n^{-1/5}C_2]$, where $0 < C_1 < 1 < C_2 < \infty$. We assume, too, that $m/n$ has a finite, nonzero limit, and that $f$ and $g$ cross at a unique point $y$ in $\mathcal{I}$, at which $\Delta'(y) \neq 0$ and the curvatures of $f$ and $g$ have different signs. The argument we shall employ to prove that $P_1$ fails in this case can be used to show that it fails more generally.

Put

$$S_0 = m^{-1} \sum_i I\{\Delta(X_i) < 0, X_i \in \mathcal{I}\} \quad \text{and} \quad U(h_1, h_2) = S(h_1, h_2) - S_0.$$

It is straightforward to show that $E\{S(h_1, h_2)\} = s(h_1, h_2) + o(\eta)$, and, of course, $S_0$ does not depend on $h_1$ and $h_2$. We shall prove that $\text{var}\{U(h_1, h_2)\}$ is asymptotic to $n^{-1}$ multiplied by a bounded function which depends nondegenerately on $(v, w) = (n^{1/5}h_1, n^{1/5}h_2)$. Call this property $P_2$, and note that $\eta^2 = o(n^{-1})$ uniformly in $h_1, h_2 \in \mathcal{H}$. It may also be proved that $U(h_1, h_2)$ is asymptotically normally distributed, and converges weakly to a Gaussian process indexed by $(v, w) \in [C_1, C_2]$. These results imply that $P_1$ fails.

Note that $\text{var}\{U(h_1, h_2)\}$, being the variance of a sum, can be expanded as a sum of diagonal terms, plus a double series in off-diagonal terms. It is relatively straightforward to show that the sum of diagonal terms equals $o(\eta)$. Therefore, it suffices to show that $P_2$ applies to the double series in off-diagonal terms contributing to the variance. That quantity equals $(1 - m^{-1})Q$, where

$$Q = \text{cov}[I\{\widehat{\Delta}_{f,-1}(X_1), X_1 \in \mathcal{I}\} - I\{\Delta(X_1) < 0, X_1 \in \mathcal{I}\},$$
$$I\{\widehat{\Delta}_{f,-2}(X_2), X_2 \in \mathcal{I}\} - I\{\Delta(X_2) < 0, X_2 \in \mathcal{I}\}],$$

and so it is adequate to prove that $P_2$ applies to $Q$.

Define $\xi = \{(n-1)h_1\}^{-1}$, $\bar{f}(x) = \xi \sum_{i \neq 1,2} K\{(x - X_i)/h_1\}$, $\delta_1(x_1, x_2) = \xi K\{(x_1 - x_2)/h_1\}$, $\delta_2(u) = \xi K(u)$, $p_j = P\{\bar{f}(x_j) - \hat{g}(x_j) + \delta_1(x_1, x_2) < 0\}$,

$q_j = P\{\hat{f}_{-j}(x_j) - \hat{g}(x_j) < 0\}$ and $r_j = I\{\Delta(x_j) < 0\}$. Let $\mathcal{K}(x_1)$ denote the set of $u$ such that $x_1 - hu \in \mathcal{I}$, and put $h = h_1$. In this notation,

$$Q = \int_\mathcal{I} \int_\mathcal{I} \{(p_1 - r_1)(p_2 - r_2) - (q_1 - r_1)(q_2 - r_2)\} f(x_1) f(x_2)\, dx_1\, dx_2$$

$$= h \int_\mathcal{I} \int_{\mathcal{K}(x_1)} \{(a_1 - a_2)(b_1 - b_2)$$
$$- (a_1 - a_3)(b_1 - b_3)\} f(x_1) f(x_1 - hu)\, dx_1\, du,$$

where $a_1 = P\{\bar{f}(x_1) - \hat{g}(x_1) < 0\} - r_1$, $a_2 = P\{-\delta_2(u) < \bar{f}(x_1) - \hat{g}(x_1) < 0\}$,

$$a_3 = P\{-\delta_1(x_1, X_2) < \bar{f}(x_1) - \hat{g}(x_1) < 0\},$$
$$b_1 = P\{\bar{f}(x_1 - hu) - \hat{g}(x_1 - hu) < 0\} - I\{\Delta(x_1 - hu) < 0\},$$
$$b_2 = P\{-\delta_2(u) < \bar{f}(x_1 - hu) - \hat{g}(x_1 - hu) < 0\},$$
$$b_3 = P\{-\delta_1(x_1 - hu, X_2) < \bar{f}(x_1 - hu) - \hat{g}(x_1 - hu) < 0\}.$$

It may thus be shown that

$$Q \sim h \int_\mathcal{I} \int_{\mathcal{K}(x_1)} \{(a_3 - a_2) b_1 + (b_3 - b_2) a_1\} f(x_1) f(x_1 - hu)\, dx_1\, du$$

$$\sim 2h \int_\mathcal{I} \int_{\mathcal{K}(x_1)} (b_3 - b_2) a_1 f(x_1)^2\, dx_1\, du$$

$$\sim -2h \int_\mathcal{I} \int_{\mathcal{K}(x_1)} a_1 b_2 f(x_1)^2\, dx_1\, du.$$

In the last-written integral, change variable from $x_1$ to $z$, where $x_1 = y + h^2 z$. Then, for arbitrarily small $\varepsilon > 0$, $Q$ is asymptotic to

$$-h^3 f(y)^2 \int_{|u| \leq n^\varepsilon} \int_{|z| \leq n^\varepsilon} P\{-\delta_2(u) < \bar{f}(y + h^2 z - hu)$$
$$- \hat{g}(y + h^2 z - hu) < 0\}$$
$$\times [I\{\Delta(y + h^2 z) > 0\}$$
$$- P\{\bar{f}(y + h^2 z) - \hat{g}(y + h^2 z) > 0\}]\, du\, dz.$$

(7.1)

The probability that occurs as a factor in the integral at (7.1) is asymptotic to $(mh)^{-1/2}$ multiplied by a nondegenerate function of $(v, w) = (n^{1/5} h_1, n^{1/5} h_2)$. The factor within square brackets in (7.1) is asymptotic to another such function. Hence, $Q$ is asymptotic to $h^3/(mh)^{1/2} \asymp n^{-1}$, multiplied by a function of $(v, w)$, as had to be proved.

**Acknowledgments.** We are grateful to two reviewers for helpful comments.

Centre for Mathematics
  and its Applications
Australian National University
Canberra, ACT 0200
Australia
e-mail: halpstat@pretty.anu.edu.au

Department of Statistics
Hankuk University
  of Foreign Studies
Mohyun, Yongin 449-791
Korea
e-mail: khkang@hufs.ac.kr